\documentclass[12pt]{article}
\usepackage[reqno]{amsmath}
\usepackage{amssymb, amsthm, enumerate}
\usepackage{graphicx}
\usepackage{tikz}
\usetikzlibrary{decorations.markings}
\usepackage{array}
\usepackage{float}
\usepackage{wrapfig}
\usepackage{url}

\usepackage{fullpage}
\expandafter\ifx\csname pdfpagebox\endcsname\relax\else
\fi

\makeatletter
\newtheorem*{rep@theorem}{\rep@title}
\newcommand{\newreptheorem}[2]{%
\newenvironment{rep#1}[1]{%
 \def\rep@title{#2 \ref{##1}}%
 \begin{rep@theorem}}%
 {\end{rep@theorem}}}
\makeatother

{

}

\newreptheorem{theorem}{Theorem}

\newcommand\cref[1]{Corollary~\ref{cor:#1}}

\def\sqr#1#2{{\vbox{\hrule height.#2pt
    \hbox{\vrule width.#2pt height#1pt \kern#1pt
        \vrule width.#2pt}\hrule height.#2pt}}}
\def\eqed{\sqr53}
\def\qed{%
    \ifmmode\eqno\eqed
    \else\nobreak\ \hfill\eqed\medbreak\fi}

\DeclareMathOperator{\tr}{tr}

\newcommand{\footremember}[2]{%
    \footnote{#2}
    \newcounter{#1}
    \setcounter{#1}{\value{footnote}}%
}
\newcommand{\footrecall}[1]{%
    \footnotemark[\value{#1}]%
} 

\title{Quantum Walks on Generalized Quadrangles}

\author{%
  Chris Godsil\footremember{uw}{Department of Combinatorics \& Optimization, University of Waterloo, Waterloo, ON N2L 3G1.  \texttt{\{cgodsil, tmyklebu\}@uwaterloo.ca} }%
  \and Krystal Guo\footnote{Department of Mathematics, Simon Fraser University, Burnaby, B.C. V5A 1S6. Currently at Department of Combinatorics \& Optimization, University of Waterloo, Waterloo, ON N2L 3G1. \texttt{kguo@uwaterloo.ca}  }
  \and Tor G.~J.~Myklebust\footrecall{uw}
  }

\begin{document}
\maketitle

\begin{abstract}
	We study the transition matrix of a quantum walk on strongly regular graphs. It is proposed by Emms, Hancock, Severini and Wilson in 2006, that the spectrum of $S^+(U^3)$, a matrix based on the amplitudes of walks in the quantum walk, distinguishes strongly regular graphs. We probabilistically compute the spectrum of the line intersection graphs of two non-isomorphic generalized quadrangles of order $(5^2,5)$ under this matrix and thus provide strongly regular counter-examples to the conjecture. 
\end{abstract}

\section{Introduction}\label{sec:intro}

A discrete-time quantum walk is a quantum process on a graph whose state vector is governed by a matrix, called the transition matrix. In \cite{ESWH, EHSW06} Emms, Severini, Wilson and Hancock propose that the quantum walk transition matrix can be used to distinguish between non-isomorphic strongly regular graphs. After experiments on a large set of graphs, no strongly regular graph was known to have a cospectral mate with respect to this invariant. 
In this paper we will compute the spectrum of $S^+(U^3)$ for two particular non-isomorphic graphs and show that they are not distinguished by  the spectrum of $S^+(U^3)$. 

A \textsl{discrete quantum walk} is a process on a graph $X$ governed by a unitary matrix, $U$, which is called the \textsl{transition matrix.} For $uv$ and $wx$ arcs in the digraph of $X$, the transition matrix is defined to be:
\[
 U_{wx,uv} = \begin{cases} \frac{2}{d(v)} &\text{ if } v=w \text{ and } u \neq x ,\\
\frac{2}{d(v)} -1  & \text{ if } v=w \text{ and } u = x, \\
0 &\text{ otherwise.} \end{cases}
\]


Let $U(X)$ and $U(H)$ be the transition matrices of quantum walks on $X$ and $H$ respectively. Given a matrix $M$, the \textsl{positive support} of $M$, denoted $S^+(M)$, is the matrix obtained from $M$ as follows:
\[ (S^+(M))_{i,j} = \begin{cases} 1 & \text{if } M_{i,j} >0\\
0 & \text{otherwise.}\end{cases}
\]

It is easy to see that if $X$ and $H$ are isomorphic  regular graphs, then $S^+(U(X)^3)$ and $S^+(U(H)^3)$ are cospectral.  For convenience, we will define $S := S^+(U(X)^3)$ and we will write $S$ or $S^+(U^3)$ to mean $S^+(U(X)^3)$ when the context is clear.  The authors of \cite{EHSW06, ESWH} propose that this spectrum is also a complete invariant for strongly regular graphs; they conjecture that the spectrum of the matrix $S^+(U(X)^3)$ distinguishes strongly regular graphs. A graph $X$ on $n$ vertices is \textsl{strongly regular} if it is neither complete nor empty, each vertex has $k$ neighbours, each pair of adjacent vertices has $a$ common neighbours and each pair of non-adjacent vertices has $c$ neighbours. The tuple $(n,k,a,c)$ is said to be the \textsl{parameter set} of $X$. 

In his Ph.D. thesis \cite{Sm12}, Jamie Smith constructs an infinite family of graphs which are not distinguishable by the procedure of Emms et al. There graphs are not strongly regular but are close, in that they have diameter two and four eigenvalues. The eigenvalues of $S^(U)$ and $S^+(U^2)$ were studied by two of the authors in \cite{GG11}. Pairs of small regular (but not strongly regular) counterexamples are given in \cite{KG10}. Hadamard graphs of order $n$ are also not distinguished by the procedure of Emms et al. for $n = 4,8,16,20$ and we conjecture that it is true for all $n$. 

In this article, we give strongly regular counterexamples to the conjecture of Emms et al. by finding a pair of non-isomorphic strongly regular graphs with parameter set $(756, 130, 4, 26)$ which have the same spectrum with respect to $S^+(U^3)$. These strongly regular graphs are the line intersection graphs of the two generalized quadrangles of order $(5^2, 5)$. The line intersection graphs of the two generalized quadrangles of order $(5^2, 5)$ are not distinguished by the procedure of Emms et al.

Since the transition matrices of these graphs are $98280\times 98280$, our computation of their minimal polynomials were done probabilistically. We then determined the eigenvalues and their multiplicities of both matrices, given the minimal polynomials. 

\section{Generalized quadrangles}\label{sec:hgrs}

The spectrum of $S^+(U^3)$ distinguishes strongly regular graphs for many parameter sets. In particular, the conjecture of Emms et al.~was checked for the small strongly regular graphs on less than or equal to $64$ vertices found in \cite{TS}. Note that the collection of strongly regular graphs in \cite{TS} is not complete for graphs up to $64$ vertices; for example, the class of strongly regular graphs with parameter set $(57, 24, 11, 9)$, consisting of graphs constructed from Steiner triple systems $S(2,3,19)$ is missing. Nevertheless, the procedure of Emms et al. distinguishes many classes of strongly regular graphs and we are motivated to search for strongly regular graphs with more regularity. 

It is known for a strongly regular graph that when the Krein bound holds with equality for the diagonal Krein parameter, the second neighbourhood of any vertex is also strongly regular. The parameter set $(756,130,4,26)$ is the smallest parameter set with a pair of non-isomorphic graphs having the property of vanishing Krein parameter. See \cite{ABsite}. We focus on strongly regular graphs with vanishing Krein parameter since the Hadamard graphs, which were also not distinguished by the procedure of Emms. et al.~are distance-regular graphs with vanishing Krein parameter.

A generalized quadrangle of order $(s,t)$ is an incidence structure where every point lies on $s+1$ lines and every line contains $t+1$ points. We construct the \textsl{line intersection graph} of a generalized quadrangle by taking the lines to be the vertices and two lines are adjacent if they have a common point. The line intersection graph of a generalized quadrangle of order $(s,t)$ is strongly regular with parameter set $((t+1)(st+1),s(t+1),t -1,s+1)$. There are two known non-isomorphic generalized quadrangles of order $(5^2,5)$ which are known in the literature as $H(3,5^2)$ and $\text{FTWKB}(5)$. 

\section{Eigenvalue computations}\label{sec:eigs}

In this section, we describe our computations for the line intersection graphs of the two generalized quadrangles of order $(5^2, 5)$, namely $H(3,5^2)$ and $\text{FTWKB}(5)$.

Let $X$ be the line intersection graph of $H(3,5^2)$ and let $Y$ be the line intersection graph of $\text{FTWKB}(5)$. Using Wiedemann's algorithm, we probabilistically compute that the minimal polynomials of both $S^+(U(X)^3)$ and $S^+(U(Y)^3)$, modulo every prime between $1999999000$ and $1999999180$, is as follows:
\begin{equation}\label{eq:minpoly}
\begin{split}
 (x-1) &(x-15)(x-125)(x-127)(x-68005)(x+25) (x+23)(x+9) \\
 & (x^2-5426x+7128229)
 (x^3+799x^2+122869x-7632765).
 \end{split}
\end{equation}
See  \cite{HaJoSa14} for analysis of Wiedemann's algorithm. Since the minimal polynomial is square-free, $S^+(U(X)^3)$ and $S^+(U(Y)^3)$ are both diagonalizable over each finite field for which the computation was done. 

We further deterministically computed $S^+(U(X)^3)^2$ and $S^+(U(Y)^3)^2$ and obtained the following:
\[ \tr(S^+(U(X)^3)) = \tr(S^+(U(Y)^3)) = 98280,\]
\[ \tr(S^+(U(X)^3)^2) = \tr(S^+(U(Y)^3)^2) = 6670853280,\]
\[\tr(S^+(U(X)^3)^3) = \tr(S^+(U(Y)^3)^3) =  318986389121400,\]
and \[\tr(S^+(U(X)^3)^4) = \tr(S^+(U(Y)^3)^4) = 21401273663621790120.\]


Let $x_i$ be the multiplicity of the $i$th factor appearing in (\ref{eq:minpoly}) as a factor of the characteristeristic polynomial of $S^+(U(X)^3)$. 
Since $S^+(U(X)^3)$ and $S^+(U(Y)^3)$ are both irreducible matrices with entries in $\{0,1\}$ the Perron-Frobenius theorem implies that the largest eigenvalue in amplitude ($68005$) has multiplicity $1$ and so $m_5 = 1$. We probabilistically compute that $m_9 = 105$ and $m_{10} = 680$; we generated $2000$ random Krylov spaces modulo $1999999151$ and   only generated $105$ eigenvectors for the $9$th factor and $680$ eigenvectors for the $10$th factor. We get the following system of linear equations
\[ { \scriptsize
\begin{split}
m_1+m_2+m_3+m_4+m_5+m_6+m_7+m_8+2m_9+3m_{10} &= 98280 \\
m_1+15m_2+125m_3+127m_4+68005m_5-25m_6-23m_7-9m_8+5426m_9-799m_{10} &= 98280 \\
m_1+225m_2+15625m_3+16129m_4+4624680025m_5& \\ 
 +625m_6+529m_7 +81m_8+15185018m_9+392663m_{10} &= 6670853280 \\
m_1+3375m_2+1953125m_3+2048383m_4+314501365100125m_5 & \\-15625m_6 -12167m_7-729m_8+43716137114m_9-192667111m_{10} &= 318986389121400 \\
m_1+50625m_2+244140625m_3+260144641m_4+21387665333634000625m_5 &\\ +390625m_6 +279841m_7+6561m_8+128961474307442m_9+99596332307m_{10} &= 21401273663621790120 \\
\end{split} }
\]
We substitute the values of $m_5$, $m_9$ and $m_{10}$ and simplify to obtain the following linear system of $5$ equations in $7$ variables:
\begin{equation}\label{eq:solve}
\begin{split}
59241m_1 +17575m_7 +72896m_8 &= 2544438125 \\
-10780m_2  +1665m_7  +4556m_8 &= 88344500 \\
 8525m_3   +570m_7  +1088m_8 &= 74452850 \\
-11172m_4   +703m_7  +1340m_8 &= -20869720 \\
12350m_6 +10545m_7  +2278m_8 &= 626911750.
\end{split}
\end{equation}

We deterministically computed lower bounds of the remaining multiplicities; we generated Krylov spaces at random and found $2000$ linearly independent eigenvectors for each remaining eigenvalue modulo $1999999151$. We obtain that $m_i \geq 2000$ for $i \in \{1,2,3,4,6,7,8\}$. Solving (\ref{eq:solve}), we find only one positive integer solution satisfying this condition:
\[
m_1 = 15625,\, m_2= 2625,\, m_3 = 4914,\, m_4 = 5460,\, m_6= 24570,\, m_7 = 27300,\, m_8 = 15625.
\]

The same computation were done for $S^+(U(Y)^3)$ and the same arguments follow and so $S^+(U(X)^3)$  and $S^+(U(Y)^3)$  are cospectral. 


\end{document}